\documentclass[10pt]{article}

\usepackage{natbib}
\usepackage{amsmath,amssymb,amsthm}

\newcommand{\qbin}[2]{\genfrac{[}{]}{0pt}{}{#1}{#2}}
\newcommand{\gp}[3]{\qbin{#1}{#2}_{#3}}

\usepackage{mathrsfs}

\newtheorem{theorem}{Theorem}

\theoremstyle{remark}
\newtheorem{remark}[theorem]{Remark}

\numberwithin{theorem}{section}

\allowdisplaybreaks

\title{On Integer partitions\\ and the\\ Wilcoxon rank-sum statistic}
\author{Andrew V. Sills}
\date{\today}

\begin{document}

\maketitle

\begin{center}
Georgia Southern University\\Department of Mathematical Sciences\\
Statesboro and Savannah, Georgia, USA
\end{center}


\begin{abstract}
In the literature, derivations of exact null distributions of rank-sum statistics is often avoided
in cases where one or more ties exist in the data.  By deriving the null distribution
in the no-ties case with the aid of classical $q$-series results of Euler and Rothe,
we demonstrate how a natural generalization of the method may be employed
to derive exact null distributions even when one or more ties are present in the
data.  It is suggested that this method could be implemented in a computer
algebra system, or even a more primitive computer language, so that the 
normal approximation need not be employed in the case of small sample sizes,
when it is less likely to be very accurate.   Several algorithms for determining exact
distributions of the rank-sum statistic (possibly with ties) have been given in the literature
(see \citet{SR86} and \citet{MBMLK16}), but none seem as simple as the procedure
discussed here which amounts to multiplying out a certain polynomial, extracting co\"efficients, and 
finally dividing by a binomal co\"efficient.
\end{abstract}

\section{Background on rank-sum tests}
\citet{W45} introduced the nonparametric rank-sum test as an alternative to
the $t$-test for difference of population means with two unpaired samples:
samples of size $n_1$ and $n_2$ are taken\footnote{\cite{W45} only
explicitly considers the case where $n_1 = n_2$.}, the data are ordered from
$1$ through $n_1+n_2$, and then the test statistic is taken to be the
sum of the ranks in one of the samples.
In the later literature, Wilcoxon's rank-sum test statistic is often denoted $W$,
and that convention will be followed here.

 \cite{W45} discusses the possibility of ties among the data: in the example in the
left column of page 80, the 12th and 13th smallest measurements are identical,
so he assigns the mean rank $12.5$ to both.  

  Many of the classic texts on nonparametric statistics avoid saying too much about
the exact distribution of $W$ in the presence of ties.  
For instance,~\citet[p. 217]{C80} advises
use of the exact rank-sum test statistic ``[i]f there are no ties, or just a few ties'', and recommends use of a normal approximation ``[i]f there are many ties.'' 
Guidelines for how few are
``few'' and how many constitute ``many'' are not provided.

\citet{MW47} extended Wilcoxon's work by allowing the two sample sizes to be
unequal, and calculating moments of their $U$-statistic, defined
in~\citet[p. 51]{MW47} as
\begin{equation} U := n_1 n_2 + \frac{n_2 (n_2 + 1)}{2} - W. 
\end{equation}  \citet[p. 51]{MW47} assume continuous data,
so the question of ties does not arise, or rather arises with a probability of $0$. 
The statistics $U$ and $W$ are ultimately equivalent; we will work with $W$ here.

\citet{L61} discusses exact and approximated distribution of the Wilcoxon rank sum
statistic with ties. \citet{G05} considers ``fuzzy'' p-values for various non-parametric
tests when ties are present in the data.

A number of statistical software programs, including R and JMP, automatically 
utilize a normal approximation for the rank-sum test when ties are present, 
even for small sample sizes.   In particular R returns a warning message that it
cannot compute an exact p-values with ties.
An algorithm for finding the exact null distribution
of the rank-sum statistic is given by~\citet{SR86} and is implemented in the
R package \texttt{coin}.   Another algorithm, based on dynamic programming, is
given by~\citet{MBMLK16}.

  The distribution of the rank-sum statistic where no ties are
  present, i.e. in sets where all data values are unequal, will be discussed here in light of partitions 
and strict-partitions, and some well-known formulas in the theory of $q$-series.
This perspective suggests how to extend to the case where ties are present
in the data, necessitating the introduction of a generalization of partitions
where half-odd integers are allowed as parts and repeated parts need to
be distinguished.  Accordingly, we will review the necessary basics of integer partition 
theory in the next section.  The necessary generalization will be motivated in context and
introduced in section~\ref{rswt} below.

\section{Integer partitions and strict-partitions}
The distribution of the Wilcoxon rank-sum statistic will be discussed as an application of
the theory of integer partitions.  In cases where ties exist in the observed data, a
generalization of the concept of integer partition must be adopted.  This generalization
follows naturally once the generating-function method of enumerating ordinary 
integer partitions is understood.

To prepare for the discussion of the rank-sum statistic, it will be useful to review
some of the basic theory of integer partitions, and to fix some notation.
\subsection{Integer partitions}
A \emph{partition} of an integer $w$ is a multiset $\lambda$ 
of positive integers whose sum
is $w$.  Each element of $\lambda$ is a \emph{part} of the partition $\lambda$,
and the number of parts of $\lambda$ is the \emph{length} of $\lambda$.  
For example, there are seven partitions of $5$:
\[ (1,1,1,1,1), (1,1,1,2), (1,1,3), (1,2,2), (1,4), (2,3), (5) . \]

 Let $p_\ell(w,\mathscr{S})$ denote the number of partitions of $w$ of length
\emph{at most} $\ell$ and with all parts in the set $\mathscr{S}$.
Thus, for example, $p_3 (5, \{1,2,3,4\}) = 4 .$

Take $\ell$ and $m$ to be fixed nonnegative integers.
It is well known that the ordinary power series generating function for
$p_\ell(w, \{ 1,2,3, \dots, m \})$ is given by
\begin{equation} \sum_{w=0}^{\ell m}  p_{\ell}(w, \{ 1,2,\dots, m \}) q^w = \gp{m+\ell}{\ell}{q}, 
\end{equation}
 where 
\begin{equation} \label{gpdef} \gp{m+\ell}{\ell}q = \prod_{i=1}^\ell 
\frac{1-q^{m+i}}{1-q^i} 
\end{equation}
is the $q$-binomial co\"efficient.  
Although not obvious from the definition~\eqref{gpdef}, $\gp {\ell+m}{\ell}q$ 
is a polynomial in $q$ of
degree $\ell m$, that satisfies recurrences analogous to the Pascal triangle
recurrence for the ordinary binomial co\"efficients: see~\citet[Chapter 3]{A76}.

\subsection{Strict-partitions}
A \emph{strict-partition} of the integer $w$ is a partition with no repeated parts,
this it is a \emph{set} of positive integers that sum to $w$.  Of the seven partitions of $5$,
three of them are strict-partitions:
\[ (1,4), (2,3), (5). \]
Let $s(w, \ell, \mathscr{S})$ denote the number of strict-partitions of $w$ with 
length exactly equal to $\ell$
and parts that lie in the set $\mathscr{S}$.  Thus $s(5, 2, \{1,2,3,4\} ) = 2$.
If we wish to emphasize the distinction be partitions and strict-partitions, we
may refer to the former as ``unrestricted partitions.''

As stated by~\citet{M16}, it turns out that
\begin{equation} \label{pe}
 s(w,\ell, \{1,2,3,\dots,m\} ) = p_{\ell} \left(w - \frac{\ell(\ell + 1)}{2},  
 \{ 1,2,\dots, m - \ell\} \right)  .
 \end{equation}
There is a natural bijection between these two classes of partitions.
\citet[p. 82]{W45} mentions~\eqref{pe}, but misstates it slightly.

\section{The rank-sum statistic $W$ and partition generating functions}
\subsection{Rank-sums where no ties occur in the data} 
In the usual nonparametric alternative to the two-sample $t$-test, we have
two samples of size $n_1$ and $n_2$ respectively. 
The data from both samples combined are ordered from smallest to largest
and assigned ranks $1, 2, 3, \dots, n_1 + n_2$.
 For one of the samples,
say the one of size $n_1$, we let the random variable $W$ equal the sum
of the ranks from this sample.  Under the null hypothesis that both samples
come from an identical distribution,
\begin{equation} P(W = w) =  \frac{ s(w, n_1, \{1, 2, 3, \dots, n_1 + n_2 \} ) }{ \binom{n_1+n_2}{n_1}},
\end{equation}
for $w = n_1(n_1 + 1)/2 ,  n_1(n+ 1)/2 + 1, \dots,  n_1(n_1 + 2n_2 + 1)/2$;  and $0$ otherwise.
To find $s(w, n_1, \{ 1,2, \dots, n_1 + n_2 \})$, Wilcoxon used a result
from~\cite{M16}, which we stated above as Eq.~\eqref{pe},
and then calculates this latter quantity via a recurrence~\cite[p. 82]{W45}.

 What we wish to point out here is that in the 1740s, 
~\citet[Chapter XVI]{E88} (and see also ~\citet[Chapter 2]{A76}) showed that
 $$s(w, n_1 , \{ 1, 2, \dots, n_1 + n_2 \})$$ is the co\"efficient of $z^{n_1} q^w$ in
 the expansion of the product
 \begin{equation} \label{Euler} 
  (1+zq)(1+zq^2)(1+zq^3) \cdots (1 + zq^{n_1 + n_2} )  = \prod_{i\in 
  \mathscr{S}} (1+zq^i),
  \end{equation} where the set $\mathscr{S} = \{ 1,2,3, \dots, n_1 + n_2 \}$,
 and that this interpretation suggests an extension that allows us to produce
the exact null distribution for $W$ when one or more ties occur in the data.

Further, in~\cite{R11}, the earliest known incarnation of the $q$-binomial theorem
is given (cf.~\cite[p. 490, Cor. 10.2.2 (c), with $x = -zq$]{AAR99}):
\begin{equation} \label{RotheQBT}
  \prod_{i=1}^N (1+zq^i) = \sum_{k=0}^N z^k \gp{N}{k}{q} q^{k(k+1)/2} ,
\end{equation}
where $N = n_1 + n_2$ is the total size of both samples together.
Thus, Rothe's formula~\eqref{RotheQBT} shows that the co\"efficient of 
$z^k$ is easily
extracted from Euler's product~\eqref{Euler}.

Eq.~\eqref{RotheQBT} is a $q$-analog of the binomial theorem because
setting $q=1$ 
results in the ordinary binomial theorem in the form
\begin{equation} (1+z)^N = \sum_{k=0}^N \binom{N}{k} z^k. \label{BT} \end{equation}
  
  Let us demonstrate the derivation of the distribution of $W$ in the case of
a small example:  suppose $n_1 = 2$, $n_2 = 3$.
 Then 
 \begin{multline}  \label{expansion23}
    (1+zq)(1+zq^2)(1+zq^3)(1+zq^4)(1+zq^5) \\ \phantom{ } \\
 =  z^5 q^{1+2+3+4+5} \\
 +z^4(q^{2+3+4+5} + q^{1+3+4+5} + q^{1+2+4+5} +q^{1+2+3+5}+ q^{1+2+3+4}) \\
 +z^3\left( q^{3+4+5} + q^{2+4+5} + q^{1+4+5} + q^{2+3+5} + q^{1+3+5} + q^{2+3+4} \right.
 \\  \left. + q^{1+3+4} +  q^{1+2+5} + q^{1+2+4} + q^{1+2+3}\right) \\
 \left(q^{4+5}+q^{3+5}+ q^{2+5} + q^{3+4} +q^{1+5} + q^{2+4}+ 
  q^{1+4} + q^{2+3} +q^{1+3}+q^{1+2} \right) z^2\\
   +\left(q^5+q^4+q^3+q^2+q\right) z + 1\\  \phantom{ } \\
 = q^{15} z^5+\left(q^{14}+q^{13}+q^{12}+q^{11}+q^{10}\right) z^4\\
 +\left(q^{12}+q^{11}+2 q^{10}+2
   q^9+2 q^8+q^7+q^6\right) z^3\\+\left(q^9+q^8+2 q^7+2 q^6+2 q^5+q^4+q^3\right) z^2\\
   +\left(q^5+q^4+q^3+q^2+q\right) z + 1,
\end{multline} 
where in the middle member, the relevant strict-partitions generated are 
displayed in the exponents of $q$.

Taking into account that $\binom 52 = 10$, and 
since the co\"efficient of $z^2$ is $q^9+q^8+2 q^7+2 q^6+2 q^5+q^4+q^3$,
we can simply read off the probability mass function (pmf) of $W$ under the null distribution:
\begin{equation} P(W=w) = \left\{  
\begin{array}{ll}
{1}/{10}, & \mbox{if $w = 3,4,8, 9$;}\\
{2}/{10}, & \mbox{if $w = 5, 6, 7;$} \\
0, & \mbox{otherwise}
\end{array}
\right. .\end{equation}

  Notice that, while we only extract the co\"efficient of $z^{n_1}$, the co\"efficients
of the powers of $q$ of $z^{n_2}$ are the same (actually, in reverse order).  
Further, in a hypothetical distributional calculation, the ``unused'' terms are 
actually not ``wasted,'' as they contain the distributional information for 
\emph{all} possible pairs $(n_1, n_2)$ for a given total of $n_1 + n_2$
observations across the two samples.  In our example, we now have the information to 
construct the distribution for the $(n_1,n_2) = (1,4)$ case as well as the $(2,3)$ case, and this 
information could be stored for future use.

By~\eqref{RotheQBT}, the pmf of $W$ in the case of sample size $n_1$
is encoded in the co\"efficients of the polynomial 
\[ q^{n_1(n_1 + 1)/2} \gp {n_1 + n_2}{n_1}q. \]
Thus the expected value $E(W)$ may be obtained, with the aid of logarithmic 
differentiation, by
\begin{align}
E(W) & = \sum_{w=n_1(n_1+1)/2}^{n_1(n_1 + 2 n_2 + 1)/2} w P(W=w) \notag \\
         & = \binom{n_1+n_2}{n_1}^{-1} \frac{d}{dq} q^{n_1 (n_1 + 1)/2} \gp{n_1 + n_2}{n_1}{q} \Bigg |_{q=1} \notag \\
         &= \frac{1}{2} n_1 (n_1 + n_2 + 1).
\end{align}

An even more tedious calculation reveals
\begin{align}
E(W^2) & = \sum_{w=n_1(n_1+1)/2}^{n_1(n_1 + 2 n_2 + 1)/2} w^2 P(W=w) \notag \\
         & = \binom{n_1+n_2}{n_1}^{-1} 
         \frac{d}{dq} \left( q\frac{d}{dq} q^{n_1 (n_1 + 1)/2} \gp{n_1 + n_2}{n_1}{q} \right) \Bigg |_{q=1} 
           \notag \\
         &= \frac{1}{12} n_1 (n_1 + n_2 + 1)\big(  n_2 + 3n_1(n_1 + n_2 + 1) \big).
\end{align}
 
 Therefore,
 \begin{equation}
    \mathrm{Var}(W) = E(W^2) - \big( E(W) \big)^2 = \frac{n_1 n_2 (n_1 + n_2 + 1)}{12}.
 \end{equation}

\subsection{Rank-sums with ties in the data: an example} \label{rswt}
Let us modify the preceding example, 
supposing now that the second and third
smallest data values are the same.  We would then assign a rank of $2.5$ to both
of these measurements.  Thus, our multiset of ranks is now $(1, 2.5, 2.5, 4, 5)$ 
and no longer the set $\{ 1,2,3,4,5 \}$.  Better yet, since our collection of ranks
contains two instances of $2.5$ of which we will need to keep track, let us label
them $2.5_1$ and $2.5_2$ and now we may consider our collection of ranks
a set once again, with two distinguished instances of $2.5$ (both of which are
equal in magnitude and therefore
contribute the same value to a rank-sum $w$):
$\{ 1, 2.5_1, 2.5_2, 4, 5 \}$.  Assuming that we still have $n_1 = 2$ and $n_2 = 3$,
let us adapt the product~\eqref{expansion23} to our current circumstance.
We have under the null hypothesis 
that $\binom 52 P(W = w)$ is the co\"efficient of $z^2 q^w$ in the
expansion of
\begin{multline} \label{GenEuler}
  (1+zq)(1+zq^{2.5_1})(1+zq^{2.5_2})(1+zq^4)(1+zq^5)\\ \phantom{ } \\
  = 1 + z(q^1 + q^{2.5_1} + q^{2.5_2} + q^4 + q^5) + \\
  (q^{1 + 2.5_1} + q^{1+2.5_2} + q^{2.5_1 + 2.5_2} + q^{1+4} + q^{1+5} +
    q^{2.5_1 + 4} + q^{2.5_2 + 4}  \\
    + q^{2.5_1 + 5} + q^{2.5_2 + 5} + q^{4+5})z^2 \\
    +(q^{1+2.5_1 + 2.5_2} + q^{1 + 2.5_1 + 4}  + q^{1 + 2.5_2 + 4} 
    + q^{1 + 2.5_1 + 5}  + q^{1 + 2.5_2 + 5} 
    + q^{2.5_1 + 2.5_2 + 4}\\  + q^{2.5_1 + 2.5_2 + 5}  + q^{1+4+5} 
    + q^{2.5_1 + 4 +5}  + q^{2.5_2 + 4 + 5}) z^3\\
    + (q^{1+2.5_1+2.5_2+4} + q^{1+2.5_1+2.5_2+5} + q^{1+2.5_1+4+5} + q^{1+2.5_2+4+5} +
     q^{2.5_1 + 2.5_2 + 4 + 5})z^4 \\
     + q^{1+2.5_1 + 2.5_2 + 4 + 5} z^5  \\  \phantom{ } \\
 = 1 + (q + 2q^{2.5} + q^4 + q^5)z + (2q^{3.5} + 2q^5 + q^6 + 2q^{6.5} 
            + 2q^{7.5} + q^9) z^2 \\
       +(q^6 + 2q^{7.5} + 2q^{8.5} + q^9 + 2q^{10} + 2q^{11.5})z^3 \\
       +(q^{10} + q^{11} + 2q^{12.5} + q^{14}) z^4 + q^{15}z^5.         
\end{multline}
Thus, the principle behind Euler's two-variable polynomial generating function
product~\eqref{Euler} remains intact, but the rank-sums encountered 
in the expansion are no
longer strict-partitions, as parts will repeat and repeated ranks must be 
distinguished from one another \emph{and} the parts are no longer necessarily
positive integers.   Let us call a sum $w$ of $n_1$ ranks with the possibility of
ties a \emph{quasipartition of $w$ with length $n_1$ from the multiset $\mathscr{R}$
consisting of the list of ranks, with tied ranks distinguished by subscripts}.  The number of such quasipartitions will be denoted
$ Q(w, n_1, \mathscr{R} ).$
In the context of the preceding example, under the null hypothesis,
\begin{equation} \label{2t} 
P(W = w) = \frac{Q(w, 2, \{ 1, 2.5_1, 2.5_2, 4, 5 \}) }{\binom{5}{2}}, 
\end{equation}
for $w \in \{ 3.5, 5, 6, 6.5, 7.5, 9 \}$; and $0$ otherwise.
Referring to expansion of the product~\eqref{GenEuler}, by reading off the
numerical co\"efficients of $q^w$ in the co\"efficient of $z^2$,~\eqref{2t} is
in fact
\begin{equation} P(W=w) = \left\{  
\begin{array}{ll}
{1}/{10}, & \mbox{if $w =6, 9$;}\\
{2}/{10}, & \mbox{if $w = 3.5, 5, 6.5, 7.5;$} \\
0, & \mbox{otherwise}
\end{array}
\right. .\end{equation}

Direct calculation shows that $E(W) = 6$ and $E(W^2) = 38.85$, so $\mathrm{Var}(W) = 2.85$.

\begin{remark}
 In the case of $m$-way ties in the data, the resulting average ranks are 
integers if $m$ is odd and half-odd-integers if $m$ is even.  Thus, no matter
how many ties are encountered, no fractions with denominator greater
than $2$ will appear among the multiset $\mathscr{R}$ of ranks.  
Accordingly, we could, in~\eqref{GenEuler} set $q = x^2$ and instead 
expand $(1+zx^2)(1+zx^5)^2(1+zx^8)(1+zx^{10})$, once again obtaining
a polynomial in two variables.  Using the natural bijection, 
\begin{equation} Q(w, 2, \{1, 2.5_1, 2.5_2, 4, 5\}) = Q(2w, 2, \{2, 5_1, 5_2, 8, 10 \}). 
\end{equation}
Any advantages of working with a true polynomial (in $z$ and $x$) over
a polynomial in $z$ and $\sqrt{q}$ must be weighed against the 
advantages of having the true rank-sums represented in the exponents
of $q$ as in the latter case.  Some experimentation in 
\emph{Mathematica} suggests
that to increase speed, it may well be advantageous to transform 
the half-integer powers of $q$ into equivalent integer powers of $x$ for the
purpose of computation, and then translate back later as necessary.

\end{remark}

\subsection{The rank-sum statistic with ties, in general}
Given the multiset $\mathscr{R}$ of ranks possibly including ties,
generalize Euler's product~\eqref{Euler} to 
\begin{equation} \mathcal{E}(z,q):= \underset{\mbox{\tiny with multiplicities}}
{\prod_{i\in\mathscr{R}}} ( 1 + zq^i ). 
\end{equation}
Then the co\"efficient of $z^{n_1}$,
\begin{equation} \frac{1}{n_1!}
\frac{\partial^{n_1}}{\partial z^{n_1}}\mathcal{E}(z,q)\Big |_{z=0}, 
\end{equation}
is a polynomial in $\sqrt{q}$ that encodes the null distribution of $W$
in that $P(W = w)$ is the co\"efficient of $q^w$ divided by $\binom{n_1+n_2}{n_1} $.
Note that this is just a way to provide a ``formula'' for extracting the 
polynomial in $\sqrt{q}$ that is the co\"efficient of $z^{n_1}$ in the
expansion of the product $\mathcal{E}(z,q)$.  It ``differentiates away'' the
powers of $z$ less than $n_1$ and removes the powers of $z$ greater
than $n_1$ by setting $z=0$.  Surely, in a computational setting, this would 
\emph{not} be an efficient way to isolate those terms in which the
power of $z$ is $n_1$.

In fact, under the null hypothesis of identical distributions, we could 
isolate the co\"efficient of $z^{n_1} q^w$ in order to write 
\begin{align*}
  P(W = w) &=  \binom{n_1+n_2}{n_1}^{-1} 
 \frac{1}{\Gamma(w+1) n_1!} \frac{\partial^w}{\partial q^w} \left[
\frac{\partial^{n_1}}{\partial z^{n_1}}
   \left( \mathcal{E}(z,q)\Big |_{z=0} \right) \right]_{q=0}\\
   & =   
 \frac{n_2!}{\Gamma(w+1) (n_1+n_2)!} \frac{\partial^w}{\partial q^w} \left[
\frac{\partial^{n_1}}{\partial z^{n_1}}
   \left( \mathcal{E}(z,q)\Big |_{z=0} \right) \right]_{q=0},
\end{align*}
where the outer derivative is a fractional derivative if $w\not\in\mathbb Z$, see, e.g.,
\cite{OS74}.

The expected value and variance of $W$ are respectively
\begin{align}
   E(W) &= \frac{n_1(n_1 + n_2 + 1)}{2}, \\
   \mathrm{Var}(W) &= \frac{n_1 n_2}{n_1+n_2 -1} \left(  \frac{1}{n_1+n_2}\sum_{i=1}^{n_1+n_2} 
   R_i^2 -
    \frac{(n_1+n_2+1)^2}{4} \right),
\end{align}
where $R_1, R_2, \dots, R_{n_1+n_2}$  are the complete list of ranks including tied ranks with 
appropriate multiplicity.  So, for the example in section~\ref{rswt}, $R_1 = 1$,
$R_2 = R_3 = 2.5$, $R_4 = 4$ and $R_5 = 5$.  

\section{Some remarks on computation}
For given sample sizes $n_1$ and $n_2$ with certain patterns of tied ranks within
the data, the example in section 2.2 demonstrates that the exact distribution 
could be easily calculated with a computer algebra system such as Maple,
Mathematica, or Sage, by extracting the appropriate co\"efficients from 
the two variable expansion of an easily formed finite product.
If a faster computation is desired, the programming could be done in, say, C or 
$\mathrm{C}{++}$,
where the polynomial in $z$ and $\sqrt{q}$ is realized as a two-dimensional array
of the relevant co\"efficients.

For the case with no ties, the computation would be even more straightforward,
thanks to~\eqref{RotheQBT}.

\section{The unconditional null distribution of $W$}
Notice that the distributional results thus far considered, are conditioned on the 
sample, in the sense that \emph{a priori} a researcher will not know whether ties
will occur in the observed data, and if they do, what the pattern of ties will be.

Continuing with the case where the total number of observations $N = 5$ and
there are $n_1 = 2$ observations
in the first group, let us enumerate all possible patterns of ties that can occur.
\vskip 5mm
\begin{tabular}{| l |c|c|}
\hline
binary &corresponding        & possible sums of $n_1 = 2$ summands \\
label    &ranks      & (with associated probability) \\
\hline \hline
0000 = 0& 3, 3, 3, 3, 3  & 6 (1) \\ \hline
0001 = 1 & 2.5, 2.5, 2.5, 2.5, 5 &  5 (0.6); 7.5 (0.4)\\ \hline
0010 = 2 & 2, 2, 2, 4.5, 4.5  &  4 (0.3); 6.5 (0.6); 9 (0.1) \\ \hline
0011 = 3 & 2, 2, 2, 4, 5 & 4 (0.3); 6 (0.3); 7 (0.3); 9 (0.1) \\  \hline
0100 = 4 &1.5, 1.5, 4, 4 , 4 & 3 (0.1); 5.5 (0.6); 8 (0.3) \\   \hline
0101 = 5 & 1.5, 1.5, 3.5, 3.5, 5 & 3 (0.1); 5 (0.4); 6.5 (0.2);\\
& &  7 (0.1); 8.5 (0.2) \\  \hline
0110 = 6 & 1.5, 1.5, 3, 4.5, 4.5 & 3 (0.1); 4.5 (0.2); 6 (0.4); \\
  & & 7.5 (0.2); 9 (0.1) \\ \hline
0111 = 7 &1.5, 1.5, 3, 4, 5 & 3 (0.1); 4.5 (0.2); 5.5 (0.2); \\
         &    & 6.5 (0.2); 7 (0.1); 8 (0.1); 
9 (0.1)\\ \hline
1000 = 8 &1, 3.5, 3.5, 3.5, 3.5 &  4.5 (0.4); 7 (0.6) \\ \hline
1001 = 9 & 1, 3, 3, 3, 5 & 4 (0.3); 6 (0.4); 8 (0.3) \\ \hline
1010 = 10 & 1, 2.5, 2.5, 4.5, 4.5 & 3.5 (0.2); 5 (0.1); 5.5 (0.2); \\
 & &  7 (0.4); 9 (0.1) \\ \hline
 1011 = 11 & 1, 2.5, 2.5, 4, 5 & 3.5 (0.2); 5 (0.2); 6 (0.1); \\
   & & 6.5 (0.2); 7.5 (0.2); 9 (0.1) \\ \hline
 1100 = 12 & 1, 2, 4, 4, 4 &   3 (0.1); 5 (0.3); 6 (0.3); 8 (0.3) \\ \hline
 1101 = 13 & 1, 2, 3.5, 3.5, 5 & 3 (0.1); 4.5 (0.2); 5.5 (0.2); \\
   & &  6 (0.1); 7 (0.2); 8.5 (0.2) \\ \hline
 1110 = 14 & 1, 2, 3, 4.5, 4.5 & 3 (0.1); 4 (0.1); 5 (0.1); 5.5 (0.2); \\ 
   &  & 6.5 (0.2); 7.5 (0.2); 9 (0.1) \\ \hline
 1111 = 15 & 1, 2, 3, 4, 5 & 3 (0.1); 4 (0.1); 5 (0.2); 6 (0.2); \\
     &  &  7 (0.2); 8 (0.1); 9 (0.1) \\ \hline
\end{tabular}
\vskip 3mm
And so there are $2^{N-1}$ possible collections of patterns of ranks, each
canonically labeled with a bit string of length $N-1$: $00000\dots0$ through $1111\dots1$.  Let us say that rank pattern $k$ occurs with probability $p_k$
for $k = 0, 1, 2, \dots, 2^{N-1}-1 $.
We must have 
$0 \leq p_k \leq 1$ for each $k$; and $$\sum_{k=0}^{2^{N-1}-1} p_k = 1.$$

In the example given in the table above, if $p_{15} = 1$ and all other 
$p_k = 0$, this is the situation where we are assuming there will be no ties.  Under this
assumption, the null distribution of $W$ is
\begin{equation}
P(W=w) = \left\{  
\begin{array}{ll}
{1}/{10}, & \mbox{if $w = 3, 4, 8, 9$;}\\
{2}/{10}, & \mbox{if $w = 5, 6, 7;$} \\
0, & \mbox{otherwise.}
\end{array}
\right. 
\end{equation}

If we go to the other extreme and assume that each of the $16$ possible 
rank patterns are equally likely, then each $p_k = 1/16$, and 
the null distribution of $W$ is as follows:
\begin{equation} P(W=w) = \left\{  
\begin{array}{ll}
{1}/{20}, & \mbox{if $w = 3,  9$;}\\
1/40, & \mbox{ if $w = 3.5, 8.5$;}\\
11/160, & \mbox{ if $w = 4, 8$;} \\
1/16, & \mbox{ if $w = 4.5, 7.5$;} \\
19/160, & \mbox{ if $w = 5, 7 $;} \\
7/80 , & \mbox{ if $w = 5.5, 6.5$;} \\
7/40,  & \mbox{ if $w = 6$;} \\
0, & \mbox{otherwise.}
\end{array}
\right. \end{equation}

Finally, if we do not specify values for $p_0, p_1, \dots, p_{15}$, the unconditional
null distribution of $W$ is given by the following pmf:
\begin{multline*} 
P(W=w)\\ = \left\{  
\begin{array}{ll}
0.1(p_4 + p_5 + p_6 + p_7+ p_{12} + p_{13}+p_{14} + p_{15}), & \mbox{if $w = 3$;}\\
0.2(p_{10} + p_{11}), & \mbox{if $w = 3.5$;}\\
0.3(p_2 + p_3+ p_9) + 0.1(p_{14} + p_{15}), & \mbox{if $w = 4$;}\\
0.2(p_6 + p_7 + 2p_8 + p_{13}), , & \mbox{if $w = 4.5$;} \\
0.6p_1 + 0.4p_5 + 0.1p_{10} + 0.2p_{11} +0.3p_{12} + 0.1p_{14} + 0.2p_{15}, & 
\mbox{if $w = 5$;} \\
0.6 p_4 + 0.2 (p_7 + p_{10} + p_{13} + p_{14}) & \mbox{if $w = 5.5$;} \\
p_0 + 0.3 (p_3 + p_{12})+ 0.4 (p_6 + p_9) + 0.1 (p_{11} + p_{13}) + 0.2p_{15}
& \mbox{if $w = 6$;} \\
0.6 p_2 + 0.2(p_5 + p_7  + p_{11} + p_{14}) , & \mbox{if $w = 6.5$;} \\
0.3 p_3 + 0.1 (p_5 + p_7 ) + 0.6 p_8 + 0.4 p_{10} + 0.2 (p_{13} + p_{15})
 & \mbox{if $w = 7$;} \\
0.4 p_1 + 0.2 (p_6 + p_{11} + p_{14}), &\mbox{if $w = 7.5$;} \\
0.3 (p_4 + p_9 + p_{12}) + 0.1 (p_7 + p_{15}) , &\mbox{if $w = 8$;} \\
0.2 (p_5 + p_{13}), &\mbox{if $w = 8.5$;} \\ 
0.1 (p_2 + p_3 + p_6 + p_7 + p_{10} + p_{11} + p_{14} + p_{15}); 
&\mbox{if $w = 9$;} \\
0, & \mbox{otherwise.}
\end{array}
\right. 
\end{multline*}

\section{Conclusion}
As noted earlier, a number of statistical software programs, including R and JMP, automatically 
utilize a normal approximation for the rank-sum test when ties are present, 
even for small sample sizes.  In a R, a warning is issued in this case that the p-value reported may not be accurate.  That is, it appears that the normal approximation is used even though it probably
is not a good idea to do so.   Additionally,~\citet{BLS00} discuss how different statistical packages
give different results for the same data in their 
respective implementations of the Wilcoxon rank-sum
test.   An algorithm for finding the exact null distribution
of the rank-sum statistic is given by~\citet{SR86} and is implemented in the
R package \texttt{coin}.   Another algorithm, based on dynamic programming, is
given by~\citet{MBMLK16}.

However, the alternative method described in section 3 above appears to be noteworthy as it involves
nothing deeper than polynomial multiplication.

\section*{Acknowledgment}
The author thanks Charles Champ for reading an early version of this manuscript, and
providing useful advice.  The author thanks the anonymous referees who provided many
useful suggestions for improving the manuscript.

\section*{Disclosure Statement}
The authors report there are no competing interests to declare.

\bibliographystyle{natbib-harv}

\end{document}